\documentclass[12pt,leqno]{amsart}%{article}
\usepackage{color}

\newtheorem{theorem}{Theorem}
\newtheorem{proposition}[theorem]{Proposition}
\newtheorem{corollary}[theorem]{Corollary}
\newtheorem{remark}{Remark}
\newtheorem{lemma}[theorem]{Lemma}
\newfont{\bb}{msbm10 at 12pt}

\def\pf{{\textit {Proof :} }}

\def\C{{\mathbb C}}
\def\R{\hbox{\bb R}}
\def\S{\hbox{\bb S}}
\def\P{\hbox{\bb P}}

\def\bS{{\mathbb S}}
\def\SM{{\mathbb{S} M}}

\newcommand{\bal}{\begin{align}}      \newcommand{\eal}{\end{align}}
\newcommand{\ba}{\begin{array}}      \newcommand{\ea}{\end{array}}
\newcommand{\bc}{\begin{center}}     \newcommand{\ec}{\end{center}}
\newcommand{\be}{\begin{enumerate}}  \newcommand{\ee}{\end{enumerate}}
\newcommand{\beQ}{\begin{eqnarray*}} \newcommand{\eeQ}{\end{eqnarray*}}
\newcommand{\bi}{\begin{itemize}}    \newcommand{\ei}{\end{itemize}}
\newcommand{\bt}{\begin{tabular}}    \newcommand{\et}{\end{tabular}}
\newcommand{\bdm}{\begin{displaymath}} \newcommand{\edm}{\end{displaymath}}

    \newcommand{\sm}{\bS\!\!\!/\,\!}
\newcommand{\D}{D\!\!\!\!/\,}
\newcommand{\nb}{\nabla\!\!\!\!/\,}
\newcommand{\mult}{\gamma\!\!\!/}
\def\qed{\hfill{q.e.d.}\smallskip\smallskip}
\begin{document}

\title[On a Shi-Tam type Inequality]{A holographic principle for the existence
of parallel spinor fields and an inequality of Shi-Tam type}

\author{Oussama Hijazi}
\address[Hijazi]{Institut {\'E}lie Cartan\\
Universit{\'e} Henri Poincar{\'e}, Nancy I\\
B.P. 239\\
54506 Vand\oe uvre-L{\`e}s-Nancy Cedex, France}
\email{hijazi@iecn.u-nancy.fr}

\author{Sebasti{\'a}n Montiel}
\address[Montiel]{Departamento de Geometr{\'\i}a y Topolog{\'\i}a\\
Universidad de Granada\\
18071 Granada \\
Spain}
\email{smontiel@goliat.ugr.es}

\begin{abstract}
Suppose that $\Sigma=\partial M$ is the $n$-dimensional boundary of a connected compact 
Riemannian spin manifold $( M,\langle\;,\;\rangle)$ with non-negative scalar curvature, and that the (inward) mean
curvature $H$ of $\Sigma$ is positive. We show that the first eigenvalue of the Dirac operator of the boundary corresponding to the conformal metric $\langle\;,\;\rangle_H=H^2\langle\;,\;\rangle$ is at least $n/2$ and equality holds if and only if there exists a parallel spinor field on $ M$. As a consequence, if $\Sigma$ admits an isometric and isospin
immersion $\phi$ with mean curvature $H_0$ as a hypersurface into another spin Riemannian manifold $ M_0$
 admitting  a parallel spinor field,  then
 \begin{equation} \label{HoloIneq}
\int_\Sigma H\,d\Sigma\le \int_\Sigma \frac{H^2_0}{H}\, d\Sigma
 \end{equation}
and equality holds if and only if both immersions have the same shape operator. In this case, $\Sigma$ has to be also connected. In the special case where $M_0=\R^{n+1}$, equality in (\ref{HoloIneq}) implies that $M$ is an Euclidean domain and $\phi$ is congruent to the embedding of $\Sigma$ in $M$ as its boundary. We also prove  that Inequality (\ref{HoloIneq}) implies the Positive Mass Theorem (PMT). Note that, using the PMT and the additional assumption that $\phi$ is a strictly convex embedding into the Euclidean space, Shi and Tam  \cite{ST1} proved the  integral inequality 
\begin{equation}\label{shi-tam-Ineq}
\int_\Sigma H\,d\Sigma\le \int_\Sigma H_0\, d\Sigma,
\end{equation}
which is stronger than (\ref{HoloIneq}) .	 
\end{abstract}



\thanks{This author was partially 
supported by a Spanish MEC-FEDER grant No. MTM2007-61775}

\date{October 8, 2011}        

\maketitle \pagenumbering{arabic}
 
\section{Introduction}

%\be

%\item Shi-Tam

%\item Results.
%\ee

In \cite{ST1}, Shi and Tam used the PMT to study  the boundary behavior 
of compact Riemannian spin manifolds with non-negative scalar curvature. More precisely, 
they proved the following:

\begin{theorem}\label{shi-tam-1}
Let $M$ be an $(n+1)$-dimensional  compact connected Riemannian spin  manifold 
with non-negative scalar curvature and mean convex boundary hypersurface $\Sigma$. 
If $\Sigma$ admits an isometric strictly convex immersion $\phi$  into $\R^{n+1}$, %given by the inclusion map of a strictly convex hypersurface $\Sigma\subset\R^{n+1}$, 
then %If $\Sigma$ admits an isometric and isospin immersion $\phi$ into the  Euclidean space $\R^{n+1}$, then 
\begin{equation}\label{shi-tam-Ineq}
\int_\Sigma H\,d\Sigma\le \int_\Sigma H_0\, d\Sigma,
\end{equation}
where $H$ is the mean curvature of $\Sigma$ as the boundary of $M$ and $H_0$ stands for the mean curvature of the immersion $\phi$ of  $\Sigma$ into $\R^{n+1}$. Equality holds if and only if $\Sigma$ is connected, $M$ is an Euclidean domain and the embedding of $\Sigma$ in $M$ and its immersion 
in $\R^{n+1}$ are congruent.
\end{theorem}   

For $n=2$, by the Weyl Embedding Theorem, the assumption that the boundary $\sigma$ embeds isometrically as a strictly convex hypersurface in $\R^3$ is equivalent to the fact the boundary $\sigma$ has positive Gauss curvature. Hence, Inequality (\ref{shi-tam-Ineq}) implies that positively curved mean convex boundaries in time-symmetric initial data sets, satisfying the dominant energy condition, have non-negative Brown-York Mass (see \cite{BY}). Note that, a generalization of this result, to subsets of general data sets, is given in \cite{LY1,LY2} and in \cite{WY} for the hyperbolic setting. \\

In this paper, we aim to prove the following main results:

\begin{theorem}\label{theorem1}
Let $( M,\langle\;,\;\rangle)$ be an $(n+1)$-dimensional compact  Riemannian  spin manifold  
with non-negative scalar curvature and with mean convex boundary hypersurface
$\Sigma$. Then, if $H$ denotes the mean curvature of $\Sigma$,
the first non-negative eigenvalue $\lambda_1(\D_H)$ of the Dirac operator corresponding to the conformal metric $\langle\;,\;\rangle_H=H^2\langle\;,\;\rangle$ satisfies
\begin{equation}\label{ineq-Main}
\lambda_1(\D_H)\ge \frac{n}{2}
\end{equation}
and equality holds if and only if there is a non trivial parallel spinor on $ M$. In this case, the  eigenspace corresponding to $\lambda_1(\D_H)=\frac{n}{2}$ consists of the restrictions to $\Sigma$ of the parallel spinor fields on
$ M$ multiplied by the function $H^{-\frac{n-1}{2}}$. Furthermore, the boundary hypersurface $\Sigma$ has to be connected.
\end{theorem}

\begin{theorem}\label{shi-tam-2} Under the same conditions as Theorem \ref{theorem1}, 
%Let $M$ be an compact spin Riemannian $(n+1)$-dimensional manifold with non-negative scalar curvature whose boundary hypersurface $\Sigma$ is mean convex. 
assume furthermore that  $\Sigma$  admits an isometric and isospin immersion into another $(n+1)$-dimensional Riemannian  spin  manifold  $M_0$ endowed with a non trivial parallel spinor field.% (for example, $M_0=\R^{n+1}$). 
Then Inequality (\ref{HoloIneq}) holds. 
%where $H$ is the mean curvature of $\Sigma$ as the boundary of $M$ and $H_0$ stands for the mean curvature of the immersion of  $\Sigma$ into $M_0$. 
Moreover, equality is achieved  if and only if  both immersions have the same shape operator. In this case, the boundary hypersurface $\Sigma$ must be connected.
\end{theorem}

\begin{theorem}\label{shi-tam-3} Under the same conditions as Theorem \ref{shi-tam-2}, assume that $M_0=\R^{n+1}$. Then
Inequality (\ref{HoloIneq}) holds and equality is achieved if and only if $\Sigma$ is connected, $M$ is an Euclidean domain and the embedding of $\Sigma$ in $M$ and its immersion 
in $\R^{n+1}$ are congruent.
\end{theorem}

\section{Preliminaries on spin manifolds} 
 
Let $(M, \langle\;,\;\rangle)$ be an $(n+1)$-dimensional Riemannian spin manifold, which we will suppose from now on to be connected, and  
denote by  ${\nabla}$ the Levi-Civita 
connection on its tangent bundle $TM$. We choose a spin structure
on $M$ and consider the corresponding spinor bundle 
$\SM$ which is a rank $2^{\left[\frac{n+1}{2}\right]}$
complex vector bundle. Denote by $\gamma$ the Clifford 
multiplication
\begin{equation}\label{Clm}
\gamma:\C\ell(M)\longrightarrow \hbox{End}(\SM)
\end{equation}
which is a fibre preserving algebra morphism. Then $\SM$ becomes a bundle of
complex left modules over the Clifford bundle $\C\ell(M)$ over the manifold 
$M$. When $n+1$ is even, the spinor bundle splits into the direct sum of the {\em positive}
and {\em negative} chiral subbundles
\begin{equation}\label{chiral}
\SM=\SM^+\oplus\SM^-,
\end{equation}
where $\SM ^{\pm}$ are defined to be the $\pm 1$-eigenspaces of 
the endomorphism $\gamma(\omega_{n+1})$, with $
\omega_{n+1}=i^{\left[\frac{n+2}{2}\right]}e_1\cdot e_2\cdots e_{n+1}$
the complex volume form. 

On the spinor bundle $\SM$, one has (see \cite{LM}) a natural Hermitian 
metric, denoted as the Riemannian metric on $M$ by $\langle\; ,\;\rangle$,  and  
the spinorial Levi-Civita  connection ${\nabla}$ acting on spinor 
fields. The Hermitian metric and ${\nabla}$
are compatible with the Clifford multiplication (\ref{Clm}) and compatible with each other. That is
\begin{eqnarray}
& X\langle\psi,\varphi\rangle = \langle{\nabla}_X\psi,
\varphi\rangle+\langle\psi,{\nabla}_X\varphi
\rangle&\label{comp1}\\  
&\langle\gamma(X)\psi,\gamma(X)\varphi\rangle =|X|^2
\langle\psi,\varphi\rangle&\label{comp2}\\
&{\nabla}_X\big(\gamma(Y)\psi\big) 
= \gamma({\nabla}_XY)\psi
+\gamma(Y){\nabla}_X\psi,&\label{comp3}
\end{eqnarray}
for any spinor fields $\psi,\varphi\in\Gamma(\SM)$ and any 
tangent vector fields
$X,Y\in \Gamma (TM)$. Since ${\nabla} \omega_{n+1}=0$, for
$(n+1)$ even, the decomposition (\ref{chiral}) becomes orthogonal and 
${\nabla}$ preserves this decomposition.

The Dirac operator ${D}$ on $\SM$ is the first order 
elliptic differential operator locally given by 
$$
{D}=\sum_{i=1}^{n+1}\gamma(e_i){\nabla}_{e_i},
$$
where $\{e_1,\dots,e_{n+1}\}$ is a local orthonormal 
frame of $TM$. When  $(n+1)$ is even, the Dirac operator interchanges
positive and negative spinor fields, that is,
\begin{equation}\label{Dchange}
{D} : \Gamma(\SM^\pm)\longmapsto \Gamma(\SM^\mp).
\end{equation}

\section{Hypersurfaces and induced structures}

 In this section, we compare the restriction $\sm \Sigma$ of the spinor bundle $\S M$ of 
a spin manifold $M$
to an orientable 
hypersurface $\Sigma$ immersed into $M$ and
its Dirac-type operator $\D$ to
the intrinsic spinor bundle $\S\Sigma$ of the induced spin structure on $\Sigma$  
and its fundamental
Dirac operator $D_\Sigma$. A fundamental case will be when the hypersurface $\Sigma$ is just the boundary $\partial M$ of a manifold $M$ with non empty boundary. These facts are in general well-known (see for example 
\cite{Bu,Tr, Ba2, BFGK} or our previous papers \cite{HMZ1, HMZ2, HMZ3, HMR,HM}). 
For completeness, we introduce the notations and the key facts.

Denote by $\nb$ the Levi-Civita connection associated with the 
induced Riemannian metric on $\Sigma$. The Gau{\ss} formula says that
\begin{equation}\label{riem-gaus}
\nb_XY={\nabla}_XY-\langle AX,Y\rangle N,
\end{equation}
where $X,Y$ are vector fields tangent to the hypersurface 
$\Sigma$, the vector field $N$ is a global unit field  
normal to 
$\Sigma$ and $A$ stands for the shape operator corresponding 
to $N$, that is,
\begin{equation}\label{shap-oper}
{\nabla}_XN=-A X,\qquad \forall X\in \Gamma(T\Sigma).
\end{equation}
We have that the restriction
\begin{equation}\label{defindu}
\sm\Sigma:=\SM_{|\Sigma}
\end{equation}
is a left module over $\C\ell(\Sigma)$ for the induced Clifford multiplication$$
\mult:\C\ell(\Sigma)\longrightarrow \hbox{End}(\sm\Sigma)$$
given by 
\begin{equation}\label{Clmind}
\mult(X)\psi=\gamma(X)\gamma(N)\psi
\end{equation}
for every $\psi\in\Gamma(\sm\Sigma)$ and $X\in\Gamma(T\Sigma)$
(note that a spinor field on the ambient manifold $M$ and its restriction to 
the hypersurface $\Sigma$ will be denoted by the same symbol).
Consider on $\sm\Sigma$ the Hermitian metric
$\langle\; ,\;\rangle$ induced from that of $\SM$. This metric immediately
satisfies the compatibility condition (\ref{comp2}) if one considers on $\Sigma$
the Riemannian metric induced from $M$ and the Clifford multiplication
$\mult$  defined in (\ref{Clmind}). Now the Gauss formula 
(\ref{riem-gaus}) implies that the spin connection
$\nb$ on $\sm\Sigma$ is given by the following spinorial Gauss formula
\begin{equation}\label{spin-gaus}
\nb_X\psi={\nabla}_X\psi-\frac{1}{2}\mult(A X)\psi
={\nabla}_X\psi-
\frac{1}{2}\gamma(A X)\gamma(N)\psi\,
\end{equation}
for every $\psi\in\Gamma(\sm\Sigma)$ and $X\in\Gamma(T\Sigma)$. Note  
that the compability conditions (\ref{comp1}),
(\ref{comp2}) and (\ref{comp3}) are satisfied by
$(\sm\Sigma,\mult,\langle\;,\;\rangle,\nb)$. 

Denote by ${\D}:\Gamma(\sm\Sigma)\rightarrow \Gamma(\sm\Sigma)$ 
the Dirac operator associated with the Dirac
bundle $\sm\Sigma$ over the hypersurface. It is a well-known fact that ${\D}$ is 
a first order elliptic differential operator which is formally 
$L^2$-selfadjoint. By  (\ref{spin-gaus}),
for any spinor field $\psi\in\Gamma(\SM)$, we have 
\begin{equation}\label{dira-extr1}
{\D}\psi=\sum_{j=1}^n\mult(e_j)\nb_{e_j}\psi
=\frac{n}{2}H\psi-\gamma(N)\sum_{j=1}^n\gamma(e_j)
{\nabla}_{e_j}\psi,
\end{equation}
where $\{e_1,\dots,e_n\}$ is a local orthonormal frame of $T\Sigma$ 
and $
H=\frac{1}{n}\hbox{trace\,}A$
is the mean curvature of $\Sigma$ corresponding to the 
orientation $N$. Using (\ref{spin-gaus}) and (\ref{shap-oper}), it is straightforward to see that
the skew-commu\-ta\-ti\-vity rule
\begin{equation}\label{D-commutes}
\D\big(\gamma(N)\psi\big)=-\gamma(N)\D\psi
\end{equation}
holds for any spinor field $\psi\in\Gamma({\sm}\Sigma)$. It is important to point out that, from 
this fact,  {\em the spectrum of $\D$ is always symmetric with respect to zero}, while this is the case
for the Dirac operator $D_\Sigma$ of the intrinsic spinor bundle {\em only when $n$ is even}. 
Indeed, in this case, we have an isomorphism of Dirac bundles$$
({\sm} \Sigma,\mult,\D)\equiv (\bS \Sigma,\gamma_\Sigma,D_\Sigma)
$$
and the decomposition ${\sm} \Sigma={\sm}\Sigma^+\oplus {\sm} \Sigma^-$, 
given by
$
{\sm} \Sigma^\pm:=\{\psi\in{\sm} \Sigma\,|\,i\gamma(N)\psi=\pm\psi\},
$ corresponds  
to the chiral decomposition 
of the spinor bundle $\bS \Sigma$. Hence $\D$ interchanges 
${\sm}\Sigma^+$ and 
${\sm}\Sigma^-$.

Yet when $n$
is odd the spectrum of $D_\Sigma$ has not to be symmetric. In fact, in this case, the spectrum of $\D$  is 
just the symmetrization
of the spectrum of $D_\Sigma$. This is why
the decomposition of $\S{M}$ into positive and 
negative chiral spinors induces an orthogonal 
and $\mult,\D$-invariant decomposition
$
{\sm}\Sigma={\sm}\Sigma_+\oplus {\sm}\Sigma_-$, with
${\sm}\Sigma_\pm:=(\S{M}^\pm)_{|\Sigma}$,
in such a way that
$$
({\sm}\Sigma_\pm,\mult,\D_{|{\sm}\Sigma_\pm}) \equiv (\bS \Sigma, \pm\gamma_\Sigma, \pm D_\Sigma).
$$ 
Moreover, $\gamma(N)$ interchanges the members of the decomposition and both two maps 
$
\gamma(N) : {\sm}\Sigma_\pm \longrightarrow {\sm}\Sigma_\mp
$
are isomorphisms.

Consequently, to study the spectrum of the induced operator $\D$ is equivalent to study the spectrum of the Dirac operator $D_\Sigma$ of the spin Riemannian structure induced on the hypersurface 
$\Sigma$.

\section{Conformal covariance of the Dirac operator} 
Consider a positive function $h$ on a Riemannian spin $n$-dimensional manifold $\Sigma$ and the 
corresponding conformal metric $\langle\;,\;\rangle^\star=h^2\langle\;,\;
\rangle$. We know that there exists (see \cite{Hit,Hij,BHMM}) a bundle isometry between the two spinor 
bundles $\sm \Sigma$ and $\sm^\star \Sigma$ corresponding to the same spin structure and
to the two conformally related metrics. For this reason, the two spinor bundles will 
be denoted by the same symbol $\sm \Sigma$. With this identification in mind, for the corresponding 
Clifford multiplications and spin connections, one has:
\begin{equation}\label{conf-change}
\mult^\star=h\mult,\qquad {\nb}^\star_X-{\nb}_X=
-\frac{1}{2h}\mult (X)\mult({\nabla}h)-\frac{1}{2h}
\langle X,{\nabla}h\rangle,
\end{equation}
for all $X\in \Gamma (T\Sigma)$. We can easily find from
(\ref{conf-change}) the relation between the two Dirac operators 
$\D$ and $\D^\star$ on $\sm \Sigma$
relative to the two conformally related  metrics on $M$. For any spinor field $\psi\in\Gamma(\sm \Sigma)$,
one has:
\begin{equation}\label{conf-boun-dira}
{\D}^\star\big( h^{-\frac{n-1}{2}}\psi\big)=h^{-\frac{n+1}{2}}{\D}\psi. 
\end{equation}
This conformal covariance of the classical Dirac operator of the spinor bundle 
was discovered by Hitchin (see \cite{Hit}).

\section{A Reilly inequality for manifolds with boundary}
Another key fact we will need is the following spinorial Reilly type inequality, valid when the manifold $M$ is compact. By integration of the well-known
Schr{\"o}dinger-Lichnerowicz formula over the compact $(n+1)$-dimensional Riemannian spin
manifold $M$  with boundary $\Sigma=\partial M$ and using a standard Schwarz inequality involving the lengths of the spin Levi-Civita connection and of the Dirac operator of $M$ (see for instance \cite{HMZ1}),  for any spinor field $\psi\in\Gamma(\bS M)$, one has 
\begin{eqnarray}\label{boun-weit-twis-ineq}
&{\displaystyle \int_\Sigma\big(\langle {\D}\psi,\psi\rangle
-\frac{nH}{2}|\psi|^2\big)\,d\Sigma\ge } \\ 
&{\displaystyle \frac{1}{4}
\int_ M {R}|\psi|^2\,d M-\frac{n}{n+1}
\int_ M |{D}\psi|^2\,d M},
\nonumber
\end{eqnarray}
where $R$ is the scalar curvature on $ M$. It is a well-known fact that equality 
occurs if and only if $\psi$ is a twistor-spinor (see \cite{BFGK} for the corresponding definition) on the bulk manifold $M$.

\section{A local boundary elliptic condition for the Dirac operator}\label{sec:gcbc}
As before, $\Sigma$ is the boundary hypersurface of an $(n+1)$-dimensional Riemannian 
spin compact manifold $M$. We define
two pointwise projection operators
$$
P_\pm:\sm\Sigma\longrightarrow \sm\Sigma
$$
on the induced Dirac bundle over the hypersurface, as follows
\begin{equation}\label{defP}
P_\pm=\frac{1}{2}\big(\hbox{Id}_{\sm\Sigma}\pm i\gamma(N)\big).
\end{equation}
Note that, as pointed out at the end of Section 3, when $n$ is even, these are nothing
but the projections onto the $\pm$-chiral subbundles $\sm\Sigma^\pm$, and when $n$ is odd, they are the projections onto the restrictions to $\Sigma$ of the chiral subbundles $\S M^\pm$ of the even-dimensional manifold $M$. It is immediate to see that $P_+$ and $P_-$ are selfadjoint and orthogonal to each other on every
$\sm\Sigma_p$, with $p\in\Sigma$.  We know  that 
these operators provide good (local) boundary conditions to solve equations
for the Dirac operator ${D}$ of $M$. In fact, it is well-known (at least when the metric of $M$ is cylindrical near $\Sigma$) that the boundary conditions $P_\pm=0$, which are sometimes called {\em MIT bag conditions} (\cite{CJJT,CJJTW,J}), satisfy the Lopatinsky-Shapiro condition for ellipticity (see \cite{Ho} and \cite[Chapter 18]{BW} for a definition, and \cite{HMR,HMZ3} for details and applications). Although the ellipticity of these boundary conditions $P_\pm=0$ is 
proved only in the cylindrical case and extensively used in the general case, we refer to \cite[Section 7, particularly Example 7.26]{BB} in order to check that this general use is correct (see also \cite{BC}). These facts can be summarized as follows.

\begin{proposition}\label{boun-ellip}
Let $ M$ be a compact Riemannian spin manifold with boundary 
 $\partial M=\Sigma$. Then the two orthogonal projection  
operators  $P_\pm$ acting on the spin bundle $\sm\Sigma$, defined in (\ref{defP}), induced on $\Sigma$ from the spin
bundle $\S M$ of a given spin structure on $M$,
provide (local) elliptic boundary conditions for the Dirac operator ${D}$ 
of $ M$.
\end{proposition}

\begin{remark}\label{remark 1}
{\rm Even though in the literature %about the subject 
one can frequently find  the assertion that
{\sl the classical Dirac operator ${D}$ does not admit elliptic boundary conditions 
in any dimensions and that, in fact, there are topological obstructions for its existence
when the dimension of $\Sigma$ is odd} \cite{BW, GLP, HMZ2,Se}, we have to point out that our  
local boundary conditions provided by the projections $P_\pm$ exist in each dimension without
any restriction.  This is due to the fact that, when $n$ is odd, the induced spin bundle $\sm \Sigma$ is not
the intrinsic spin bundle $\S\Sigma$ of the induced spin structure, but a direct sum of two copies of it, since 
we are considering the total spin bundle $\S M$ on the ambient manifold instead of one of its
chiral subbundles. 
Then, in this case, the Dirac operator $D$ is really a pair of classical chiral Dirac operators and
the usual restrictions do not apply to this situation. Note also that the Green integral formula
\begin{equation}\label{adjoint}
\int_ M\langle{D}\psi,\varphi\rangle\,d M-
\int_ M\langle\psi,{D}\varphi\rangle\,d M=
\int_\Sigma\langle\psi,\gamma(N)\varphi\rangle\,d\Sigma,
\end{equation}
where $\psi,\varphi\in\Gamma(\S M)$, shows that none of the conditions provided by $P_\pm$ makes ${D}$
a formally selfadjoint operator. Instead, one can easily see that the boundary {\em
realizations} $(D,P_+)$ and $(D,P_-)$ of $D$ are adjoint to each other.}
\end{remark}

The ellipticity of the boundary conditions given by $P_+$ and $P_-$ and that of the Dirac operator $D$ of $M$ guarantee 
that we may solve boundary value problems for $D$ on $M$ prescribing on the boundary $\Sigma$ the corresponding
$P_\pm$-projections of the solutions. For completeness,  we give a proof.

\begin{proposition}\label{boun-prob}
The following two types of inhomogeneous problems  for the Dirac operator ${D}$ 
of a compact Riemannian spin  manifold $M$, with 
boundary a hypersurface $\Sigma$, 
\begin{equation}\label{conf-loca-boun-cond}
\left\{
\begin{array}{rll}
{D}\psi&=\Psi \qquad& \hbox{ {\rm on} } M \\
P_\pm(\psi_{|\Sigma})&=0 \qquad& \hbox{ {\rm on} }\Sigma
\end{array}
\right. 
\end{equation}
have a unique smooth solution for any $\Psi\in\Gamma(\bS M)$.
\end{proposition}

\noindent
\pf  The two realizations of ${D}$ associated
with the two boundary conditions $P_\pm$ are the two unbounded operators
\beQ
{D}_\pm:\hbox{\rm Dom}\,{D}_\pm=\{\psi\in H^1(\bS M)\,|\,
P_\pm(\psi_{|\Sigma})=0\}\longrightarrow L^2(\bS M)
\eeQ
where $H^1(\S M)$ stands for the Sobolev space of $L^2$-spinors with weak
$L^2$ covariant derivatives (recall that such spinors have a well
defined $L^2$ trace on $\Sigma$). 
From (\ref{adjoint}) and the end of Remark \ref{remark 1}, it follows that for the adjoint, one has $(D_\pm)^* = D_\mp$. Moreover, if $\psi\in \hbox{Dom}\,{D}_\pm$ is a solution
to the corresponding homogeneous problem, that is, if it belongs to 
$\ker {D}_\pm$,
the ellipticy of both the Dirac operator $D$ and the boundary condition
$P_\pm=0$ imply regularity results from which (see \cite[Corollary 7.18]{BB}) 
one gets that $\psi$ is smooth.
On the other hand, taking $\varphi=i\psi$ in (\ref{adjoint})  and recalling that 
the metric on $\bS M$ is Hermitian, we have 
$$
0=2\int_ M\langle{D}\psi,i\psi\rangle\,d M=\int_\Sigma\langle\psi,
i\gamma(N)\psi\rangle\,d\Sigma=\mp\int_\Sigma|\psi|^2\,d\Sigma.
$$
Then one sees that the smooth harmonic spinor $\psi$ on
the compact manifold $ M$ has a vanishing trace $\psi_{|\Sigma}$ along the boundary hypersurface $\Sigma$.  But according
to \cite{Ba1}, in an $(n+1)$-dimensional manifold like $M$, the Haussdorf measure of the zero set of a smooth 
non trivial harmonic 
spinor must be less than or equal
to $n-1$. So, the spinor field $\psi$ vanishes on the whole  
of $ M$. Then 
$$
\ker{D}_\pm=\{0\}\qquad \hbox{and}\qquad \hbox{\rm coker}\,
{D}_\pm\cong\ker (D_\pm)^*=\ker{D}_\mp=\{0\}.
$$ 
Then the two realizations ${D}_\pm$ are invertible operators, hence 
if $\Psi
\in\Gamma(\bS M)$ is a smooth spinor field on $ M$, there exists a unique solution $\psi\in H^1(\bS M)$ of (\ref{conf-loca-boun-cond}). 
Now, Proposition \ref{boun-ellip} and the regularity results proved in \cite[Theorem 7.17]{BB} 
imply (cf. also \cite[Chapter 19]{BW}) the required smoothness of the solution $\psi$.

\qed

\begin{remark}\label{remark 2}
{\rm 
When the dimension $(n+1)$ of the manifold $M$ is odd, 
the final considerations made in Section 3 along with Proposition \ref{boun-prob} 
give the existence, uniqueness and the regularity to the two problems
$$
\left\{
\begin{array}{rll}
{D}\psi&=\Psi \qquad&\hbox{{\rm on}}\; M  \\
\psi_{|\Sigma}&\in \Gamma (\sm\Sigma^\mp)  \qquad&\hbox{{\rm on }} \Sigma 
\end{array}
\right. 
$$
In contrast, when  $(n+1)$ is even, we may decompose the given spinor fields $\psi$ and 
$\Psi$ according to the chiral subbundles (\ref{chiral}) of $ M$. Thus Proposition
\ref{boun-prob} solves the following boundary first order system
$$
\left\{
\begin{array}{rll}
{D}\psi_\pm&=\Psi_\mp \qquad&\hbox{{\rm on }}  M  \\
i\gamma(N)({\psi_\pm}_{|\Sigma})&=\mp {\psi_\pm}_{|\Sigma} \qquad&\hbox{{\rm on }} \Sigma, 
\end{array}
\right. 
$$
where now all the involved fields have a fixed chirality.}
\end{remark}

\begin{proposition}  \cite{HMZ3}\label{boun-prob2}
Let $ M$ be a compact Riemannian spin manifold with boundary a hypersurface
$\Sigma$. If $\varphi\in\Gamma(\sm\Sigma)$ is a smooth spinor field in the induced 
Dirac bundle and $\Psi\in\Gamma(\bS M)$, then the following boundary problem
for the Dirac operator 
$$ 
\left\{
\begin{array}{rll}
{D}\psi&=\Psi \qquad&\hbox{ {\rm on} } M   \\
P_\pm(\psi_{|\Sigma})&=P_\pm\varphi \qquad&\hbox{ \rm on }\Sigma
\end{array}
\right. 
$$
has a unique smooth solution $\psi\in\Gamma(\bS M)$.
\end{proposition}

\noindent\pf Extend $\varphi$ to a spinor field $\widehat\varphi\in\Gamma(\bS M)$ and put $\widehat\psi=\psi-\widehat\varphi$. Then solve
$$
\left\{
\begin{array}{rll}
{D}\,\widehat\psi&=-{D}\,\widehat\varphi +\Psi\qquad&\hbox{ {\rm on} } M  \\
P_\pm(\widehat\psi_{|\Sigma})&=0 \qquad&\hbox{ {\rm on} }\Sigma 
\end{array}
\right. 
$$
using Proposition \ref{boun-prob}.

\qed

\begin{remark}\label{remark 3}
{\rm 
Note that when one considers the well-known elliptic global APS boundary condition,
introduced by Atiyah, Patodi and Singer in order to study a version of the Index Theorem in the case 
of manifolds with non empty boundary (see \cite{APS}), the boundary problems
for the Dirac operator  corresponding to those solved in Propositions \ref{boun-prob} and \ref{boun-prob2}
do not necessarily have solutions. From %One may look at the works related to 
the spinorial proofs
of the Positive Mass Theorem (see for instance \cite{He,Wi}), it appears that  for solving boundary problems 
with the APS boundary condition,
it is necessary to impose some non negativity condition on the scalar curvature of $M$ and some
lower estimate on the mean curvature of the boundary hypersurface (see also \cite{HMR} for other types of boundary 
conditions). 
}
\end{remark}

\section{A holographic principle for the existence of parallel spinors}
It is by now a known approach (see \cite{HMZ2, HMZ3}) to make use of the Reilly type inequality 
(\ref{boun-weit-twis-ineq}) for a compact Riemannian spin manifold $ M$
with non-negative scalar curvature $R$, together with the solution of an appropriate boundary problem
for the Dirac operator $D$ of $M$, in order to establish a certain integral inequality 
((\ref{final-ineq}) in this case) for
the Dirac operator $\D$ of the boundary hypersurface $\partial M
=\Sigma$. Moreover, we will also assume that the inward mean curvature $H$ of $\Sigma$ is positive,  
that is, $\Sigma$ is {\em mean convex},  and so Inequality (\ref{final-ineq}) will be translated into some results about the first eigenvalue of the
Dirac operator $\D_H$, associated with the conformal metric $\langle\;,\;\rangle_H=H^2\langle\;,\;\rangle$
on $\Sigma$.
First, we need to recall the following fact: 

\begin{lemma}\label{lemma1} 
For any smooth spinor field $\psi\in\Gamma(\sm\Sigma)$ we have
$$
\int_\Sigma\langle{\D}\psi,\psi\rangle\,d\Sigma=2 \int_\Sigma
\langle{\D}P_+\psi,P_-\psi\rangle\,d\Sigma.
$$
\end{lemma}
\noindent\pf We have the pointwise orthogonal decomposition $\psi=P_+\psi+P_-\psi$.
Moreover, from  (\ref{D-commutes}) and (\ref{defP}), one 
immediately shows that
\begin{equation}\label{antiDP}
{\D}P_\pm=P_\mp{\D}.
\end{equation}
Hence, since $P_+$ and $P_-$ are orthogonal to each other,
$$
\langle{\D}\psi,\psi\rangle=\langle{\D}P_+\psi,P_-\psi\rangle +
\langle {\D}P_-\psi , P_+\psi\rangle.
$$
We conclude by noting that $\Sigma$ is compact and the operator ${\D}$ is 
formally $L^2$-selfadjoint.

\qed

\begin{proposition}\label{proposition4} Let $ M$ be a compact spin Riemannian manifold 
with non-negative scalar curvature, whose boundary hypersurface $\Sigma=\partial M$ has positive
(inward) mean curvature $H$. Let $\varphi\in\Gamma(\sm\Sigma)$ be any spinor field on the
restricted Dirac bundle. Then
\begin{equation}\label{ineq+}
0\le \int_\Sigma \big(\frac{1}{H}|\D P_+\varphi|^2-\frac{n^2}{4}H|P_+\varphi|^2\big)\,d\Sigma.
%0\le\int_\Sigma \left(\frac{|\D P_+\varphi|^2}{H}-\frac{n^2}{4} H|P_+\varphi|^2\right)\,d\Sigma.
\end{equation}
Moreover, equality holds if and only if there exists a parallel spinor field $\psi
\in \Gamma({\mathbb S} M)$ such that $P_+\psi=P_+\varphi$
along the boundary hypersurface $\Sigma$.
\end{proposition}
\noindent\pf
Unless otherwise stated,  we shall use the same symbols to indicate spinor fields defined on $M$ and their restrictions  to the boundary hypersurface $\Sigma$. Take any spinor field 
$\varphi\in\Gamma(\sm\Sigma)$ of the 
induced spin bundle on the hypersurface and consider the following boundary problem
%\begin{equation}\label{boun-prob-hom}
$$
\left\{
\begin{array}{rll}
{D}\psi&=0 \qquad&\hbox{ {\rm on} } M \\
P_+\psi&=P_+ \varphi \qquad&\hbox{ {\rm on} }\Sigma 
\end{array}
\right.
$$ 
%\end{equation}
for the Dirac operator ${D}$ and the
boundary condition $P_+$. Proposition
\ref{boun-prob2} asserts that this problem has a unique smooth solution $\psi\in
\Gamma(\bS M)$.   Putting it in the Reilly 
Inequality (\ref{boun-weit-twis-ineq}) and taking into account that we are assuming
 $R\ge 0$ on $ M$, we obtain the following key
inequality 
\begin{equation}\label{final-ineq}
0\le \int_\Sigma\big(\langle {\D}\psi,\psi\rangle
-\frac{n}{2}H|\psi|^2\big)\,d\Sigma,
\end{equation}
where, if equality is achieved, then $\psi$ is a twistor-spinor field. But it is also
harmonic, hence it is parallel.
Using Inequality (\ref{final-ineq}) combined with Lemma \ref{lemma1} above and the fact that the decomposition$$
\psi=P_+\psi+P_-\psi$$
is pointwise orthogonal, we get
\begin{equation}\label{final-ineq-P}
0\le \int_\Sigma \left(2\langle{\D}P_+\psi,P_-\psi\rangle
-\frac{n}{2}H|P_+\psi|^2-
\frac{n}{2}H|P_-\psi|^2\right)\,d\Sigma.
\end{equation}
The mean curvature $H$ being assumed positive, we consider  the obvious pointwise inequality 
\begin{eqnarray}
&{\displaystyle
0\le \big| \frac{1}{\sqrt{\frac{n}{2}H}}\D P_+\psi-\sqrt{\frac{n}{2}H} P_-\psi\big|^2
= } \nonumber \\
&{\displaystyle \frac{1}{\frac{n}{2}H}|\D P_+\psi|^2+\frac{n}{2}H|P_-\psi|^2-2\langle\D P_+\psi,
P_-\psi\rangle.}\nonumber
\end{eqnarray}
In other words, we have
\beQ
2\langle\D P_+\psi , P_-\psi\rangle -\frac{n}{2}H|P_-\psi|^2 \le \frac{1}{\frac{n}{2}H}|\D P_+\psi|^2,
\eeQ
which, when combined with Inequality 
(\ref{final-ineq-P}), implies Inequality (\ref{ineq+}).
Now, in order to study the equality case, recall that the harmonic spinor $\psi$ on $ M$ is such that %satisfied the boundary condition 
$P_+\psi=P_+\varphi$. If equality holds, we already know from
(\ref{final-ineq}) that the spinor field $\psi$ must be parallel. 

Conversely, assume that there is a parallel spinor field $\psi$ on $ M$. Then from (\ref{dira-extr1}), along the boundary $\Sigma$, we have $$
\D\psi=\frac{n}{2}H\psi\,.$$
Using the relations (\ref{antiDP}), the previous  equality splits into$$
\D P_+\psi=\frac{n}{2}HP_-\psi\qquad \text{\rm and} \qquad \D P_-\psi=\frac{n}{2}HP_+\psi.$$
From these two relations and the formal $L^2$-selfadjointness of $\D$, it is straightforward to see
that equality holds in (\ref{ineq+}) for $P_+\psi\in \Gamma(\sm\Sigma)$.

\qed

 With this, we are ready to state  the following key result: 
\begin{proposition}\label{proposition1}
Let $ M$ be a compact spin Riemannian $(n+1)$-dimensional manifold 
with non-negative scalar curvature, whose  boundary hypersurface
$\Sigma$ has positive (inward) mean curvature $H$ (that is, $\Sigma$ is {\em mean convex}). Then 
\begin{equation}\label{ineq-total}
0\le \int_\Sigma \big(\frac{1}{H}|\D \varphi|^2-\frac{n^2}{4}H|\varphi|^2\big)\,d\Sigma,
\end{equation}
for any spinor field $\varphi\in\Gamma(\sm\Sigma)$. Equality holds if and only if 
there exist two parallel spinor fields $\psi,\Psi\in\Gamma({\mathbb S} M)$ such that $P_+\psi=P_+\varphi$ and $P_-\Psi=P_-\varphi$
on the boundary.  
\end{proposition}
\noindent
\pf
Since there is an obvious symmetry between the two boundary conditions $P_+$ and $P_-$ for the Dirac operator on $ M$
(see Proposition \ref{boun-prob2}), one can repeat the proof of Proposition \ref{proposition4}
 to get the inequality corresponding  to (\ref{ineq+}) where the {\em positive} projection $P_+$ is replaced by the {\em negative} one $P_-$. Hence,
 for any spinor field $\varphi\in\Gamma(\sm\Sigma)$, we also have
\begin{equation}\label{ineq-}
0\le \int_\Sigma \big(\frac{1}{H}|\D P_-\varphi|^2-\frac{n^2}{4}H|P_-\varphi|^2\big)\,d\Sigma.
\end{equation}
Taking into account the relation (\ref{antiDP}) and the pointwise orthogonality of
the projections $P_\pm$, the sum of the two inequalities (\ref{ineq+}) and (\ref{ineq-}) yields (\ref{ineq-total}). The equality case
is a consequence of Proposition \ref{proposition4}.

\qed

\begin{remark}
{\rm Note that, in the case of equality, we cannot conclude that the two parallel spinors in Proposition \ref{proposition1} coincide. In fact, assume that
the spin manifold $M$ admits a space of parallel spinor fields with dimension at least $2$. Take two different parallel spinor fields  $\psi,\Psi\in\Gamma({\mathbb S} M)$ and define $\varphi=P_+(\psi_{|\Sigma}) + P_-(\Psi_{|\Sigma})$. Such a spinor field  on the boundary $\Sigma$ achieves the equality in Inequality (\ref{ineq-}). } 
\end{remark}

\noindent {\textit {Proof of Theorem \ref{theorem1} :} }
We consider on the boundary hypersurface $\Sigma$ the conformally modified Riemannian metric
$\langle\;,\;\rangle_H=H^2\langle\;,\;\rangle$. Using the conformal covariance of the Dirac operator  (see (\ref{conf-boun-dira})) we have that, for any
spinor field $\varphi\in\Gamma(\sm\Sigma)$,
\begin{equation}\label{eigenstar}
\D_H\varphi_H=H^{-\frac{n+1}{2}}\D\varphi,
\end{equation}
where $\varphi_H=H^{-\frac{n-1}{2}} \varphi $. Then, as the Riemannian measures of the two conformally related metrics satisfy $$
d\Sigma_H=H^nd\Sigma,$$
we obtain the equalities$$
|\D_H\varphi_H|^2d\Sigma_H=\frac{1}{H}|\D\varphi|^2d\Sigma\qquad \text{\rm and}\qquad
|\varphi_H|^2d\Sigma_H=H|\varphi|^2d\Sigma.$$
Now, it suffices to put this information into Inequality  (\ref{ineq-total}) to get
\begin{equation}\label{ineq-total*}
0\le \int_\Sigma \big(|\D_H \varphi_H|^2-\frac{n^2}{4}|\varphi_H|^2\big)
\,d\Sigma_H,
\end{equation}
which is valid for all $\varphi_H=H^{-\frac{n-1}{2}} \varphi $, with  $\varphi\in
\Gamma(\sm\Sigma)$ arbitrary. This is equivalent to the inequality $$
\lambda_k(\D_H)^2\ge \frac{n^2}{4}$$ 
for all the eigenvalues $\lambda_k(\D_H)$ of $\D$, $k\in{\mathbb Z}$. This proves Inequality (\ref{ineq-Main})  for
$\lambda_1(\D_H)$.

If equality holds in (\ref{ineq-Main}), then there is a non trivial $\varphi_H=H^{-\frac{n-1}{2}}\varphi\in\Gamma(\sm\Sigma)$ such that $\D_H\varphi_H=\frac{n}{2}\varphi_H$. From (\ref{eigenstar}), this is equivalent
to $\D\varphi=\frac{n}{2}H\varphi$. Then, it is clear that $\varphi$ satisfies the equality in (\ref{ineq-total}). Thus there exist 
two parallel spinor fields $\psi$ and $\Psi$ on $ M$
with $P_+\psi=P_+\varphi$ and $P_-\Psi=P_-\varphi$. From (\ref{dira-extr1}) and (\ref{antiDP}), it follows $$
\frac{n}{2}HP_-\psi=\D P_+\psi=\D P_+\varphi=\frac{n}{2}HP_-\varphi.$$
Hence, $P_\pm\psi=P_\pm\varphi$, and so $\varphi$ is the restriction to $\Sigma$ of the parallel spinor field $\psi$ on $ M$.

As for the connectedness of $\Sigma$ in the equality case, take a non trivial eigenspinor $\varphi_H\in\Gamma(\sm\Sigma)$ of $\D$ associated with the eigenvalue $\frac{n}{2}$. Choose now a connected component $\Sigma_0$ of $\Sigma$ where $\varphi_H$ is non trivial and define a new spinor field $\widetilde{\varphi}_H$ on $\Sigma$ in the following way:
$$
\widetilde{\varphi}_H=\left\{
\begin{array}{ll}
\varphi_H \qquad&\hbox{ {\rm on} } \Sigma_0 \\
0  \qquad&\hbox{ {\rm on} }\Sigma-\Sigma_0. 
\end{array}
\right.
$$
It is clear that $\widetilde{\varphi}_H\in\Gamma(\sm\Sigma)$ and that
$$
\D\widetilde{\varphi}_H=\frac{n}{2}\widetilde{\varphi}_H,
$$ 
that is, $\widetilde{\varphi}_H$ is another non trivial eigenspinor of $\D$ associated with the eigenvalue $\lambda_1(\D)=\frac{n}{2}$. Hence, it is the restriction to $\Sigma$ of a non trivial parallel spinor field $\psi$ on $M$ multiplied by the positive function $H^{-\frac{n-1}{2}}$. Since $M$ is connected, the length of $\psi$ is constant and so $\widetilde{\varphi}_H$ has no zeros.  This implies $\Sigma=\Sigma_0$ and 
then $\Sigma$ is connected.         

\qed

\begin{remark}
{\rm Since the Euclidean space $M=\R^{n+1}$ is a spin manifold admitting  a space of parallel spinor fields with maximal
dimension $2^{[\frac{n+1}{2}]}$, we have $\lambda_1(\D_H)=\frac{n}{2}$ for any compact mean convex {\em embedded} hypersurface $\Sigma\subset\R^{n+1}$ 
(which always bounds a compact domain), and the corresponding associated eigenspace is $2^{[\frac{n+1}{2}]}$-dimensional. Moreover, we deduce that: {\em a compact mean convex hypersurface embedded in $\R^{n+1}$ must be connected}. When the
hypersurface $\Sigma$ is allowed to have self-intersections, that is, when $\Sigma$ is an orientable hypersurface {\em immersed} in $\R^{n+1}$ with nowhere vanishing mean curvature, we  proved in \cite{HM} that only the inequality $\lambda_1(\D_H)\le \frac{n}{2}$ occurs, and the equality implies that the associated eigenspace comes from the parallel spinor fields of the Euclidean space $\R^{n+1}$ as well. } 
\end{remark}

\noindent {\textit {Proof of Theorem \ref{shi-tam-2} :} }
Choose a connected component $\Sigma_0$ of $\Sigma$ and define $\varphi\in\Gamma(\sm\Sigma)$ to be the restriction to $\Sigma_0$ of a non trivial 
parallel spinor field on $M_0$ and to be identically zero on 
$\Sigma-\Sigma_0$. Then, using (\ref{dira-extr1}), we have$$
\D\varphi=\frac{n}{2}H_0 \varphi$$
on the whole of $\Sigma$. Since $\varphi$ has a non zero constant length on $\Sigma_0$, it is sufficient  to apply Inequality (\ref{ineq-total}) to  $\varphi$ in order to get  Inequality (\ref{HoloIneq}) on the component $\Sigma_0$. Thus, the same inequality must hold on the whole of the boundary $\Sigma$. Suppose now that equality holds. By Proposition \ref{proposition1}, there exist two parallel spinor fields $\psi,\Psi\in\Gamma({\mathbb S} M)$ such that $P_+\psi=P_+\varphi$ and $P_-\Psi=P_-\varphi$. Then, using (\ref{dira-extr1}), (\ref{antiDP}) and the equality above, we have
\begin{equation}\label{eq1}
H_0P_+\varphi=\frac{2}{n}\D P_-\varphi=\frac{2}{n}\D P_-\Psi=
HP_+\Psi.   
\end{equation}
Similarly, we obtain 
\begin{equation}\label{eq2}
H_0P_-\varphi=\frac{2}{n}\D P_+\varphi=\frac{2}{n}\D P_+\psi=
HP_-\psi.
\end{equation}
Applying the operator $\D$ to the first and last terms of  (\ref{eq1}), we get $$
\mult(\nabla H_0)P_+\varphi+\frac{n}{2}H_0^2P_-\varphi = \mult(\nabla H)P_+\Psi+\frac{n}{2}H^2P_-\Psi,  
$$
and using again the equalities above, we have finally $$
\mult(\nabla H_0)P_+\varphi+\frac{n}{2}H_0^2P_-\varphi= \frac{H_0}{H}\mult(\nabla H)P_+\varphi+\frac{n}{2}H^2P_-\varphi. $$
The same argument applied to (\ref{eq2}), yields$$ 
\mult(\nabla H_0)P_-\varphi+\frac{n}{2}H_0^2P_+\varphi= \frac{H_0}{H}\mult(\nabla H)P_-\varphi+\frac{n}{2}H^2P_+\varphi. $$
The sum of the last two  formulae, implies$$
\mult(\nabla H_0)\varphi+\frac{n}{2}H_0^2\varphi= \frac{H_0}{H}\mult(\nabla H)\varphi+\frac{n}{2}H^2\varphi. $$
Since the spinor fields $\mult (\nabla H)\varphi$ and $\mult (\nabla H_0)\varphi$ are both orthogonal to $\varphi$, 
and the spinor $\varphi$ has non trivial constant length on $\Sigma_0$, we finally obtain
\beQ
\left\{
\begin{array}{rl}
H_0^2 &= H^2   \\
\nabla H_0 &= \frac{H_0}{H}\nabla H 
\end{array}
\right. 
\eeQ
on that component $\Sigma_0$. From this we conclude that $H_0$ has no zeros and so we may assume that$$
H_0=H.$$
Coming back now to (\ref{eq1}) and (\ref{eq2}), we deduce that
\beQ
\left\{
\begin{array}{rl}
P_+\varphi &= P_+\Psi  \\
P_-\varphi &= P_-\psi 
\end{array}
\right. 
\eeQ
Thus the two parallel spinor fields $\psi,\Psi\in\Gamma({\mathbb S} M)$ satisfy$$
\psi_{|\Sigma_0}=\Psi_{|\Sigma_0}=\varphi_{|\Sigma_0}.$$

By definition, on $\Sigma-\Sigma_0$ we have $\varphi_{|\Sigma-\Sigma_0}=0$. Thus$$
P_+\psi=P_+\varphi=0,\qquad P_-\Psi=P_-\varphi=0.$$
 Applying the induced Dirac operator $\D$ to these equalities and taking again into account (\ref{dira-extr1}) and (\ref{antiDP}), we have $$
0=\D P_+\psi=\frac{n}{2}HP_-\psi,\qquad 0=\D P_-\Psi=\frac{n}{2}HP_+\Psi.$$ 
We conclude that both parallel spinor fields $\psi$ and $\Psi$ on $M$ vanish on the complement of $\Sigma_0$. Since $M$ is connected, $\psi$ and $\Psi$ must have constant length
on the whole of $M$. However, this length vanishes on $\Sigma-\Sigma_0$ and is a non zero constant on $\Sigma_0$. This proves that the boundary hypersurface $\Sigma$ is connected.  

As another conclusion, we have that, if equality holds in Inequality (\ref{HoloIneq}), the mean curvatures $H$ and $H_0$ coincide and that each restriction $\varphi\in\Gamma(\sm\Sigma)$ to $\Sigma$ of a  
parallel spinor field on $M_0$ is the restriction to the boundary hypersurface $\Sigma$ of a parallel spinor field $\psi$ defined on the whole of $M$.

Now we can apply to such a $\varphi\in\Gamma(\sm\Sigma)$ the spinorial Gau{\ss} formula (\ref{spin-gaus}) for  the embedding of $\Sigma$ in $ M$ as its boundary  
 and so get$$
 \nb_X\varphi=\nabla^ M_X\psi-\frac{1}{2}\mult(A X)\varphi=-\frac{1}{2}\mult(A X)\varphi,
 $$
where $X\in\Gamma(T\Sigma)$. Using again (\ref{spin-gaus}) for the immersion of $\Sigma$ in the manifold $M_0$, we have$$  
\nb_X\varphi=\nabla^{M_0}_X\varphi-\frac{1}{2}\mult(A_0 X)\varphi=-\frac{1}{2}\mult(A_0 X)\varphi,
 $$
where  $A_0$ is the shape operator of this immersion. Since the spinor field $\varphi$ has constant length, we conclude that the two shape operators $A$ and $A_0$, corresponding to the embedding of $\Sigma$ in $ M$ and to the immersion of $\Sigma$ in $M_0$,  coincide. 

The converse is clear. If the two shape operators $A$ and $A_0$, associated to the immersions of $\Sigma$ in $M$ as its boundary and in
$M_0$ respectively, coincide, then the corresponding traces $nH$ and $nH_0$ taken with respect to the common induced metric should be equal. From here equality in
(\ref{HoloIneq}) is straightforward.

\qed

We have explicitely stated that Inequality (\ref{HoloIneq}) in Theorem \ref{shi-tam-2} above is valid when the mean convex boundary hypersurface $\Sigma$ of the Riemannian spin $(n+1)$-dimensional manifold $M$ with non-negative scalar curvature can be immersed, in an isometric and isospin way, in the Euclidean space
$\R^{n+1}$. In fact, in this particular situation, we may considerably improve the equality case. \\

\noindent {\textit {Proof  of Theorem \ref{shi-tam-3} :}}
It suffices to study the  equality case. From Theorem \ref{shi-tam-2} above, it is clear that the boundary hypersurface $\Sigma$ has to be connected. On the other hand, since the Euclidean space $\R^{n+1}$ admits a maximal number $2^{[\frac{n+1}{2}]}$ of linearly independent parallel spinor fields, we can repeat the initial argument in the proof of Theorem \ref{shi-tam-2} for each one of the restrictions to $\Sigma$ of these spinor fields. In this way, we obtain this same number $2^{[\frac{n+1}{2}]}$ of independent parallel spinor fields defined on the bulk manifold $M$. But, according to \cite{Wa}, this maximal number is only attained only by flat manifolds. Thus, $M$ is a flat manifold and so we may see it as an open set of an Euclidean quotient, that is, $M\subset\R^{n+1}/\Gamma$, where $\Gamma$ is a group of Euclidean motions acting properly and discontinuously on $\R^{n+1}$. Let $\widetilde{M}\subset\R^{n+1}$ be any connected component of the lifting of $M$ to the Euclidean space. Then any connected component of the boundary $\widetilde{\Sigma}$ of $\widetilde{M}$ is a connected hypersurface embedded in $\R^{n+1}$ covering the original hypersurface $\Sigma$. Denote by $\pi:\widetilde{\Sigma}\rightarrow \Sigma$ the corresponding projection. Then the composition $\phi\circ\pi:\widetilde{\Sigma}
\rightarrow \R^{n+1}$ is an immersion whose first and second fundamental forms coincide with those of the embedding $\iota:\widetilde{\Sigma}\subset \R^{n+1}$ as the boundary of the domain 
$\widetilde{M}$. Then we may apply the fundamental theorem of the local theory of surfaces (see, for example, \cite[Theorem 7.7, p.\! \!209]{MR}) to deduce that there exists a rigid motion
$F$ of the Euclidean space $\R^{n+1}$ such that $\phi\circ\pi=F\circ\iota$. Hence the covering map $\pi$ is injective. This means that $\widetilde{\Sigma}=\Sigma$, $\pi={\rm Id}_{|\Sigma}$ and $\phi=
F\circ\iota$.

\qed 

\begin{remark}
{\rm  It is clear that Theorem \ref{shi-tam-3}  provides an integral inequality involving the mean curvatures of two isometric (and isospin) immersions as hypersurfaces of a compact Riemannian manifold $\Sigma$ in two different ambient spaces $M$ and $\R^{n+1}$.
An inequality of this type was first obtained by Shi and Tam \cite[Theorem 4.1]{ST1} under stronger hypotheses  than that of Theorem \ref{shi-tam-3}. In fact, they
assumed that the isometric immersion $\phi$ is the inclusion map of a strictly convex hypersurface $\Sigma\subset\R^{n+1}$ (note that, due to its convexity, each connected component of $\Sigma$ must be diffeomorphic to an $n$-dimensional sphere and so $\Sigma$ admits a unique spin structure). Under this assumption, Shi and Tam proved the following inequality
\begin{equation*}%\label{shi-tam}
\int_\Sigma H\,d\Sigma\le \int_\Sigma {H_0}\, d\Sigma,
\end{equation*}
and equality holds under the same conditions as in Theorem \ref{shi-tam-3}. Observe that, combining a Schwarz inequality and 
Inequality (\ref{shi-tam-Ineq}), we get
\begin{equation}\label{chain}
\big(\int_\Sigma H_0\,d\Sigma \;\big)^2\le \int_\Sigma \frac{H^2_0}{H}\, d\Sigma \int_\Sigma H\, d\Sigma
\le \int_\Sigma \frac{H^2_0}{H}\, d\Sigma \int_\Sigma H_0\, d\Sigma.
\end{equation} 
Since $\phi$ is a strictly convex embedding, we have that the mean curvature $H_0$ of $\phi$ is positive. Thus, from the 
inequality above we deduce
\begin{equation*}%\label{hijazi-montiel-2}
\int_\Sigma H_0\,d\Sigma\le  \int_\Sigma \frac{H^2_0}{H}\, d\Sigma.
\end{equation*}
(Note that we can also deduce  Inequality (\ref{HoloIneq}) by combining the first inequality in (\ref{chain}) and Inequality (\ref{HoloIneq}), provided that the integral of $H_0$ on $\Sigma$ is positive). Inequality (\ref{HoloIneq}) along with (\ref{shi-tam-Ineq}) means that {\em the Shi and Tam inequality implies  Theorem \ref{shi-tam-3}. Note that, 
%Yet it is worth pointing out that, by way of contrast, 
Theorem \ref{shi-tam-3} is valid in a  more general setup  than that of  Inequality (\ref{shi-tam-Ineq})
(and its generalizations, such that obtained in \cite{EMW}), since no convexity assumptions on the immersion $\phi$ is
imposed.} } 
\end{remark}

\begin{remark}
{\rm   Another important remark about the relationship between  Theorem \ref{shi-tam-3} and the Shi and Tam result \cite[Theorem 4.1]{ST1}
is that a key ingredient in their proof is a version of the PMT (see \cite{SY}) for $C^2$-metrics. Furthermore,
Shi and Tam showed \cite[Theorem 5.1]{ST1} that their inequality also implies the PMT (at least in the 3-dimensional case, the
more significant from the physical point of view). Instead, {\em the proof of  Inequality  (\ref{HoloIneq}) makes no use of the PMT. However, even though it is weaker than  Inequality (\ref{shi-tam-Ineq}), it  implies the PMT,} as we will see in the following:} 
\end{remark}

\begin{corollary}\label{PMT}
For any compact connected Riemannian 3-dimensional manifold $M$ with non-negative scalar curvature and
mean convex boundary surface $\Sigma$, assume that 
$$
\int_\Sigma H\,d\Sigma\le \int_\Sigma \frac{H^2_0}{H}\, d\Sigma,
$$
where $H$ is the mean curvature of $\Sigma$ as the boundary of $M$ and $H_0$ stands for the mean curvature of any immersion $\phi$ of  $\Sigma$ into $\R^3$. Let $P$ be a complete non compact 3-dimensional Riemannian asymptotically flat manifold with finitely many ends, with non-negative integrable scalar curvature. Then the ADM mass of each end of $P$ is non-negative. 
\end{corollary}   

\noindent\pf
Without loss of generality, we may assume that $P$ has only one end, denote it by $E$. Then, we may assume that $E$ is the exterior of an Euclidean ball of $\R ^3$ and, from the work of Schoen and Yau (see the proof of \cite[Theorem 5.1]{ST1} and references therein), the assumption of asymptotically flatness may be translated into the fact that the Riemannian metric $\langle\;,\;\rangle_P$ on $E$ differs from the Euclidean metric $\langle\;,\;\rangle$ in this way:
$$ 
\langle\;,\;\rangle_P(x)=\big(1+\frac{m}{|x|}\big)\langle (x)\;,\;\rangle+h_x,\qquad\forall x\in E\subset\R^3,$$
where $h_x$ is a symmetric bilinear form satisfying
$$
\big|\frac{\partial^k h}{\partial x_i}(x)\big|=O\big(\frac1{|x|^{2+k}}\big),\qquad i=1,\dots,3,
\qquad k=0,\dots,4,$$
for all $x\in E$, and where the constant $m\in\R$ is just the ADM mass of the end $E$. Thus, we aim  to prove that $m\ge 0$.

For each $r>0$ large enough take the Euclidean sphere $\S^2_r\subset E$ of radius 
$r$ centered at the origin. Consider now the compact 3-dimensional manifold $M$ obtained by taking off from $P$ the exterior domain determined by $\S^2_r$. Then $M$ is a compact connected Riemannian manifold with non-negative scalar curvature and whose boundary is just the surface $\S^2_r$. Recall also that $M$, like all  3-dimensional manifolds, is spin. In order to apply to this manifold $M$ the integral inequality that we  assumed, we need to compute the inward mean curvature $H$ of $\S^2_r$ with respect to the metric $\langle\;,\;\rangle_N$.
In fact, this more or less straightforward computation can be seen in \cite{ST1}. More precisely, they proved that
$$
H=\frac1{r}-\frac{2m}{r^2}+O\big(\frac1{r^3}\big),\qquad K=\frac1{r^2}-\frac{2m}{r^3}+O\big(\frac1{r^4}\big),$$
where $K$ is the Gau{\ss} curvature of $\S^2_r$ with the metric induced from $\langle\;,\;\rangle_N$. This proves that, for $r$ large enough, 
not only the boundary surface $\S^2_r$ is mean convex, that is, that $H>0$, but also that it is strictly convex, that is, $K>0$. This last property
guarantees, by applying the Weyl Embedding Theorem (see \cite{We}) that the sphere $\S^2_r$ can be embedded in $\R^3$ in such a way that the metrics induced on $\S^2_r$ from $\R^3$ and the one induced from $M$ coincide. Moreover, this embedding is unique up to Euclidean congruences. This means that the mean curvature $H_0$ of this embedding is completely determined by this metric. In fact, by using the estimates for $H_0$ in terms of $K$ obtained by Weyl, Shi and Tam got the following estimate at infinity (see \cite{ST1} having in mind that their mean curvatures are not normalized) 

$$
H_0=\frac1{r}-\frac{m}{r^2}+O\big(\frac1{r^3}\big).$$
From the asymptotic expressions for $H$ and $H_0$ above, it is immediate to conclude that$$
\frac{H^2_0}{H}=\frac1{r}+O\big(\frac1{r^3}\big).$$
Finally, in order to compute the integrals on $\S^2_r$ of the two functions $H$ and ${H_0^2}/{H}$, we need some information about the Riemannian measure $d\,\S^2_r$. It is clear that the map $F:\S^2\rightarrow\S^2_r$ given by $F(y)=
ry$ for $y\in\S^2$ is a diffeomorphism %between the unit sphere $\S^2$ and  $\S^2_r$ 
and can be used as a parametrization of $\S^2_r$. We can see in \cite[(5.5)]{ST1} that$$
F_r^*\left(d\,\S^2_r\right)=\left(r^2+2mr+O(1)\right)d\,\S^2,$$ 
where $d\,\S^2$ is the Riemannian measure of the Euclidean unit sphere. Now, we may write asymptotic expressions for the integrals
on the sphere $\S^2_r$ of the three functions $H$, $H_0$ and $H^2_0/H$. In fact, we have
$$
\begin{array}{rl}
{\displaystyle \int_{{\mathbb S}^2_r}H\,d\,\S^2_r} & =  4\pi\big(r+O\big(\frac1{r}\big)\big)  \\\
{\displaystyle \int_{{\mathbb S}^2_r}H_0\,d\,\S^2_r} & =  4\pi\big(r+m+O\big(\frac1{r}\big)\big) \\\
{\displaystyle \int_{{\mathbb S}^2_r}\frac{H^2_0}{H}\,d\,\S^2_r} & =  4\pi\big(r+2m+O\big(\frac1{r}\big)\big),
\end{array}
$$
for all $r>0$ large enough. Now, one can see that  Inequality (\ref{HoloIneq}),
assumed to be true as a hypothesis, with $\Sigma=\S^2_r$, $r>0$, $r\rightarrow\infty$, implies that $m\ge 0$.
 
\qed

\begin{remark}
{\rm Note that, in Corollary \ref{PMT}, we can substitute Inequality (\ref{HoloIneq}) either by the Shi and Tam Inequality 
(\ref{shi-tam-Ineq}), which is valid only under convexity assumptions for $\Sigma$, or by Inequality (\ref{HoloIneq}), which is valid   
only when the integral on $\Sigma$ of the mean curvature $H_0$ is positive. The PMT can be deduced from any of these three inequalities for the boundary of a compact connected 3-dimensional manifold with non-negative scalar curvature. The crucial difference between them is that the Shi and Tam inequality was proved using as one of its key ingredients this Positive Mass Theorem, and also that the realm of application of the other two is wider.}   
\end{remark}
  
From Theorem \ref{shi-tam-3}, it is not difficult to deduce as well a congruence result for Euclidean immersions of mean convex boundaries of compact connected Riemannian spin manifolds with non-negative scalar curvature. This result (see  Corollary \ref{Ros-congruence}) generalizes that of Ros (\cite{R}) where the bulk manifold is supposed to have non-negative Ricci curvature. In \cite{R}, this
congruence result is presented as a generalization of an old rigidity theorem proved by Schur for plane Euclidean curves. This same generalization is obtained in \cite{HW}, where it is presented as a solution to a conjecture by Schroeder and Strake, and also in \cite[Lemma 5]{EMW}.      
   
 \begin{corollary}\label{Ros-congruence}
 Let $M$ be a compact connected spin Riemannian $(n+1)$-dimensional manifold 
with non-negative scalar curvature and mean convex boundary hypersurface $\Sigma$. Suppose that $\Sigma$ admits an isometric and isospin immersion $\phi$ into the Euclidean space $\R^{n+1}$ and that the mean curvature $H$ of $\Sigma$ as the boundary of $M$ and  the mean curvature $H_0$ of the immersion $\phi$ of  $\Sigma$ in $\R^{n+1}$ satisfy the pointwise inequality $|H_0|\le H$.
Then, $M$ is an Euclidean domain with connected boundary. Moreover, the embedding of $\Sigma$ in $M$ as its boundary and the immersion 
of $\Sigma$ in $\R^{n+1}$ are congruent. 
\end{corollary}   

\noindent\pf
Since we assume that the inequality $|H_0|\le H$ holds, we have that$$
\int_\Sigma \frac{H^2_0}{H}\, d\Sigma\le \int_\Sigma H\, d\Sigma,$$
and so equality is attained in Inequality (\ref{HoloIneq}).% Theorem \ref{shi-tam-3}. 
This finishes the proof.

\qed

\end{document}